\documentclass[a4paper,12pt]{article}
\usepackage{amsmath, amsfonts, amsthm, amssymb, mathrsfs, pb-diagram, graphics, graphicx}
\usepackage{verbatim, color, pgf, tikz, mathrsfs, fullpage}
\usepackage[numbers]{natbib}
\usetikzlibrary{patterns}
\usetikzlibrary{arrows}

\theoremstyle{definition}

\newtheorem{thm}{Theorem}

\newtheorem{que}{Question}

\newcommand{\R}{\mathbb{R}}

\DeclareMathOperator{\Aut}{Aut}
\DeclareMathOperator{\Aff}{Aff}
\DeclareMathOperator{\Proj}{Proj}
\DeclareMathOperator{\Sim}{Sim}
\DeclareMathOperator{\Euc}{Euc}

\usepackage{hyperref}
\hypersetup{colorlinks,% 
citecolor=red,% 
linkcolor=blue,% 
}
\usepackage{color}

\begin{document}
	\title{A projective analogue \\ of Napoleon's and Varignon's theorems}
	\author{Quang-Nhat Le}
	\maketitle
	
	Let us consider the following classical results.
	
	\begin{thm}[Napoleon's theorem]
		Given an arbitrary triangle. The centers of the equilateral triangles that are erected on the edges and lie outside the triangle form an equilateral triangle.
	\end{thm}
	
	\begin{thm}[Varignon's theorem]
		The midpoints of the edges of any quadrilateral form a parallelogram.
	\end{thm}
	
	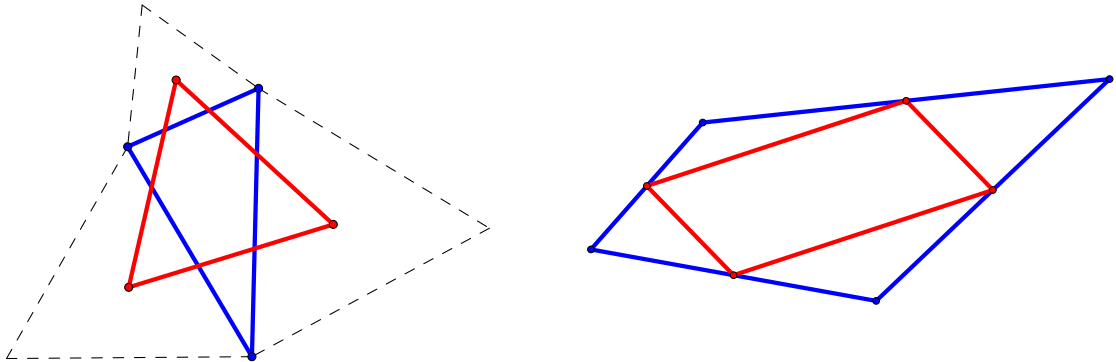
\begin{figure}[!ht]
		\centering
		
		\begin{minipage}{.5\textwidth}
			\centering
			\definecolor{ffqqqq}{rgb}{1.,0.,0.}
			\definecolor{qqqqff}{rgb}{0.,0.,1.}
			\begin{tikzpicture}[line cap=round,line join=round,>=triangle 45,x=1.0cm,y=1.0cm,scale=0.6]
			\clip(-3.266130187769224,-5.319241689249532) rectangle (10.01542794167257,5.0147579649660665);
			\draw [line width=1.6pt,color=qqqqff] (3.6094419765868118,2.1637672817774947)-- (0.7383028276882347,0.8738352003882811);
			\draw [line width=1.6pt,color=qqqqff] (0.7383028276882347,0.8738352003882811)-- (3.4638044835267388,-3.765759221382599);
			\draw [line width=1.6pt,color=qqqqff] (3.6094419765868118,2.1637672817774947)-- (3.4638044835267388,-3.765759221382599);
			\draw [dash pattern=on 4pt off 4pt] (1.0567584504979286,4.005280681829088)-- (0.7383028276882347,0.8738352003882811);
			\draw [dash pattern=on 4pt off 4pt] (1.0567584504979286,4.005280681829088)-- (3.6094419765868118,2.1637672817774947);
			\draw [dash pattern=on 4pt off 4pt] (3.6094419765868118,2.1637672817774947)-- (8.671743814206526,-0.9271217385360568);
			\draw [dash pattern=on 4pt off 4pt] (8.671743814206526,-0.9271217385360568)-- (3.4638044835267388,-3.765759221382599);
			\draw [dash pattern=on 4pt off 4pt] (3.4638044835267388,-3.765759221382599)-- (-1.9169529769026699,-3.8063156825098545);
			\draw [dash pattern=on 4pt off 4pt] (-1.9169529769026699,-3.8063156825098545)-- (0.7383028276882347,0.8738352003882811);
			\draw [line width=1.6pt,color=ffqqqq] (1.8015010849243247,2.3476277213316217)-- (5.2483300914400255,-0.8430378927137199);
			\draw [line width=1.6pt,color=ffqqqq] (5.2483300914400255,-0.8430378927137199)-- (0.7617181114374344,-2.2327465678347247);
			\draw [line width=1.6pt,color=ffqqqq] (0.7617181114374344,-2.2327465678347247)-- (1.8015010849243247,2.3476277213316217);
			\begin{scriptsize}
			\draw [fill=qqqqff] (3.6094419765868118,2.1637672817774947) circle (2.5pt);
			\draw [fill=qqqqff] (0.7383028276882347,0.8738352003882811) circle (2.5pt);
			\draw [fill=qqqqff] (3.4638044835267388,-3.765759221382599) circle (2.5pt);
			\draw [fill=ffqqqq] (1.8015010849243247,2.3476277213316217) circle (2.5pt);
			\draw [fill=ffqqqq] (5.2483300914400255,-0.8430378927137199) circle (2.5pt);
			\draw [fill=ffqqqq] (0.7617181114374344,-2.2327465678347247) circle (2.5pt);
			\end{scriptsize}
			\end{tikzpicture}
		\end{minipage}%
		\begin{minipage}{.5\textwidth}
			\centering
			\definecolor{ffqqqq}{rgb}{1.,0.,0.}
			\definecolor{qqqqff}{rgb}{0.,0.,1.}
			\begin{tikzpicture}[line cap=round,line join=round,>=triangle 45,x=1.0cm,y=1.0cm,scale=0.5]
			\clip(-2.85071246077018,-4.731408258815843) rectangle (13.087548713793064,3.143967850968349);
			\draw [line width=1.6pt,color=qqqqff] (-1.7729718627628128,-2.2862237741157676)-- (1.1683714084814552,1.078672928187679);
			\draw [line width=1.6pt,color=qqqqff] (1.1683714084814552,1.078672928187679)-- (11.87486091581059,2.2316794905154334);
			\draw [line width=1.6pt,color=qqqqff] (11.87486091581059,2.2316794905154334)-- (5.733336165452559,-3.6510070519731093);
			\draw [line width=1.6pt,color=qqqqff] (5.733336165452559,-3.6510070519731093)-- (-1.7729718627628128,-2.2862237741157676);
			\draw [line width=1.6pt,color=ffqqqq] (-0.30230022714067883,-0.6037754229640443)-- (6.521616162146023,1.655176209351556);
			\draw [line width=1.6pt,color=ffqqqq] (6.521616162146023,1.655176209351556)-- (8.804098540631575,-0.709663780728838);
			\draw [line width=1.6pt,color=ffqqqq] (8.804098540631575,-0.709663780728838)-- (1.980182151344873,-2.9686154130444384);
			\draw [line width=1.6pt,color=ffqqqq] (1.980182151344873,-2.9686154130444384)-- (-0.30230022714067883,-0.6037754229640443);
			\begin{scriptsize}
			\draw [fill=qqqqff] (-1.7729718627628128,-2.2862237741157676) circle (2.5pt);
			\draw [fill=qqqqff] (1.1683714084814552,1.078672928187679) circle (2.5pt);
			\draw [fill=qqqqff] (11.87486091581059,2.2316794905154334) circle (2.5pt);
			\draw [fill=qqqqff] (5.733336165452559,-3.6510070519731093) circle (2.5pt);
			\draw [fill=ffqqqq] (-0.30230022714067883,-0.6037754229640443) circle (2.5pt);
			\draw [fill=ffqqqq] (6.521616162146023,1.655176209351556) circle (2.5pt);
			\draw [fill=ffqqqq] (8.804098540631575,-0.709663780728838) circle (2.5pt);
			\draw [fill=ffqqqq] (1.980182151344873,-2.9686154130444384) circle (2.5pt);
			\end{scriptsize}
			\end{tikzpicture}
		\end{minipage}
		\caption{Napoleon's theorem (left) and Varignon's theorem (right).} \label{fig:NapoVari}
	\end{figure}
	
	What do these theorems have in common? We claim that they both are examples of \textit{immediately regularizing natural polygon iterations} in different planar geometries and present a new and analogous result in projective geometry.
	
	\section{Introduction}
	Polygons have always been an object of fascination of mathematicians since the ancient times. They are appealing because of their simplicity, elegance and ubiquity in mathematics. The study of geometric constructions on polygons is probably as old as mathematics itself. If such a construction creates a new polygon, we can consider it as a self-map on the space of polygons and study the dynamical behavior of its iterates. We use the term \textbf{polygon iteration} to emphasize our focus on the dynamical aspect of that construction.
	
	One classical example of polygon iterations is the midpoint map $\mathcal{M}$ which constructs from a polygon $P$ a new polygon $Q$ whose vertices are the midpoints of the edges of $P$. The map $\mathcal{M}$ is affinely natural, i.e. it commutes with all affine transformations, and is known to have connections to finite Fourier analysis \cite{Terras1999} and outer billiards \cite{Tabachnikov1999, Troubetzkoy2009}. The most interesting dynamical aspect of the map $\mathcal{M}$ is perhaps the fact that, under the iteration of $\mathcal{M}$, the affine shape of $P$ will converge to that of a regular polygon, for almost any starting polygon $P$. We say that the iterates of $\mathcal{M}$ \textbf{regularize} $P$. 
	
	\begin{figure}[!ht]
		\centering
		\definecolor{ffqqqq}{rgb}{1.,0.,0.}
		\definecolor{qqqqff}{rgb}{0.,0.,1.}
		\begin{tikzpicture}[line cap=round,line join=round,>=triangle 45,x=1.0cm,y=1.0cm,scale = 0.7]
		\clip(-3.5776859504132235,-2.519504132231412) rectangle (6.009090909090915,5.29867768595042);
		\draw [line width=1.6pt,color=qqqqff] (-2.44,-0.02)-- (-2.24,3.34);
		\draw [line width=1.6pt,color=qqqqff] (-2.24,3.34)-- (0.46,5.02);
		\draw [line width=1.6pt,color=qqqqff] (0.46,5.02)-- (3.88,3.22);
		\draw [line width=1.6pt,color=qqqqff] (3.88,3.22)-- (5.18,-0.86);
		\draw [line width=1.6pt,color=qqqqff] (5.18,-0.86)-- (0.2,-2.32);
		\draw [line width=1.6pt,color=qqqqff] (-2.44,-0.02)-- (0.2,-2.32);
		\draw [line width=1.6pt,color=ffqqqq] (-2.34,1.66)-- (-0.89,4.18);
		\draw [line width=1.6pt,color=ffqqqq] (-0.89,4.18)-- (2.17,4.12);
		\draw [line width=1.6pt,color=ffqqqq] (2.17,4.12)-- (4.53,1.18);
		\draw [line width=1.6pt,color=ffqqqq] (4.53,1.18)-- (2.69,-1.59);
		\draw [line width=1.6pt,color=ffqqqq] (2.69,-1.59)-- (-1.12,-1.17);
		\draw [line width=1.6pt,color=ffqqqq] (-1.12,-1.17)-- (-2.34,1.66);
		\begin{scriptsize}
		\draw [fill=qqqqff] (-2.44,-0.02) circle (2.5pt);
		\draw [fill=qqqqff] (-2.24,3.34) circle (2.5pt);
		\draw [fill=qqqqff] (0.46,5.02) circle (2.5pt);
		\draw [fill=qqqqff] (3.88,3.22) circle (2.5pt);
		\draw [fill=qqqqff] (5.18,-0.86) circle (2.5pt);
		\draw [fill=qqqqff] (0.2,-2.32) circle (2.5pt);
		\draw [fill=ffqqqq] (-2.34,1.66) circle (2.5pt);
		\draw [fill=ffqqqq] (-0.89,4.18) circle (2.5pt);
		\draw [fill=ffqqqq] (2.17,4.12) circle (2.5pt);
		\draw [fill=ffqqqq] (4.53,1.18) circle (2.5pt);
		\draw [fill=ffqqqq] (2.69,-1.59) circle (2.5pt);
		\draw [fill=ffqqqq] (-1.12,-1.17) circle (2.5pt);
		\end{scriptsize}
		\end{tikzpicture}
		\caption{The midpoint map $\mathcal{M}$ on a hexagon. ($P$ is blue and its image $Q$ is red.)} \label{fig:Midpoint}
	\end{figure}
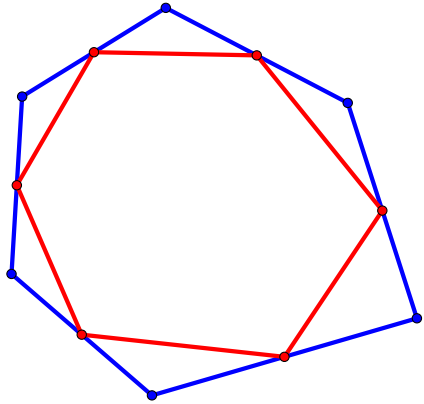
	
	The above statement is vacuous when $P$ is a triangle simply because all triangles are affinely equivalent to each other. In the first nontrivial case of quadrilaterals, Varignon's theorem asserts that the midpoint map $\mathcal{M}$ sends any quadrilateral to a parallelogram, which is affinely regular because it can be mapped to a square by an affine transformation. In other words, the map $\mathcal{M}$ \textbf{immediately regularizes} quadrilaterals. Therefore, we can interpret Varignon's theorem as saying that the midpoint map is an affinely natural and immediately regularizing polygon iteration on quadrilaterals.
	
	Analogously, we can see Napoleon's theorem as a statement in similarity geometry describing a natural and immediately regularizing polygon iteration on triangles. 
	
	How can we extend these results to projective geometry? First, let us explain what immediately regularizing natural polygon iterations are in projective geometry.
	
	\section{Klein geometries}
	
	In geometry, Felix Klein's Erlangen program \cite{Klein1893} is one of those ``grand schemes of things'' that propose a radically new perspective that puts various geometric theories in the same tent and connects them to big contemporary ideas. At the heart of the Erlangen program is the interpretation of a geometry $\mathbf{G}$ on a space $X$ as the action of the Lie group $\Aut(\mathbf{G})$ of automorphisms of $\mathbf{G}$ on $X$. For example, according to Klein, the Euclidean geometry of the plane $\R^2$ should be thought of in terms of the action of the group $\Euc(\R^2)$ of Euclidean transformations on $\R^2$. Therefore, the study of geometries can be reduced to that of the corresponding Lie groups and their actions, which was a burgeoning field in Klein's time and is now a classical and well-studied area of mathematics. 
	
	Among the Klein geometries of the plane, there is a tower of supergeometries of Euclidean geometry $\textbf{Euc}^2$:
	\begin{equation*} \label{eq:tower}
		\textbf{Sim}^2 \subset \textbf{Aff}^2 \subset \textbf{Proj}^2,
	\end{equation*}  
	which comprises of similarity geometry, affine geometry and projective geometry. Their automorphism groups can be ordered by inclusion as following:
	\[ \Sim_2(\R) \subset \Aff_2(\R) \subset \Proj_2(\R). \]
	Strictly speaking, the Lie group $\Proj_2(\R)$ acts on the projective plane $\R P^2$, which is an extension of the plan $\R^2$, but this technical point will not noticeably affect our discussions.
	
	\subsection{Natural constructions}
	
	Klein's view of geometries as group actions provides a notion of natural geometric constructions. A geometric construction, particularly a polygon iteration, is called \textbf{natural} if it commutes with any automorphism of the underlying Klein geometry. Therefore, Napoleon's construction is natural in similarity geometry and so is Varignon's in affine geometry.
	
	Natural polygon iterations can also be considered as self-maps on the \textbf{moduli space} of polygons, which is defined as the quotient of the space of polygons by the diagonal action of the automorphism group $\Aut(\mathbf{G})$.
	
	\subsection{Regular polygons}
	
	In a planar Klein geometry, a polygon is called \textbf{regular} if it can be carried by an automorphism to a regular polygon, defined in the Euclidean sense. For example, parallelograms are affinely regular. Furthermore, a polygon iteration on $n$-gons, i.e. polygons with $n$ vertices, is \textbf{immediately regularizing} if it transforms an arbitrary $n$-gon to a regular $n$-gon.
	
	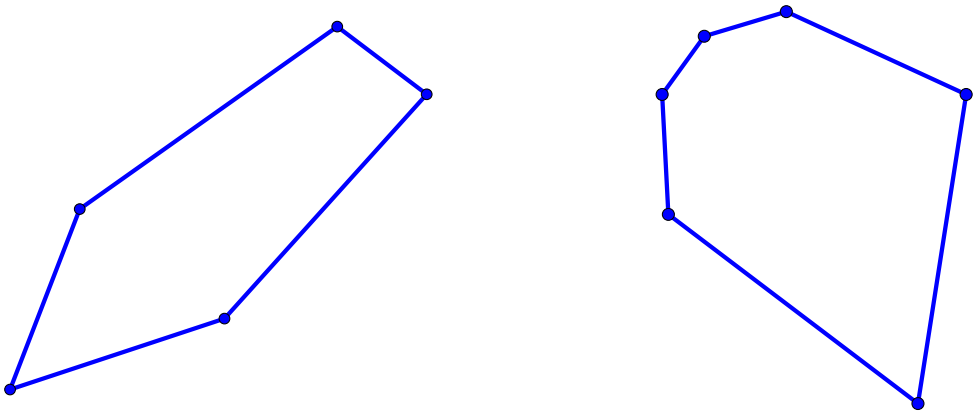
\begin{figure}[!ht]
		\centering
		
		\begin{minipage}{.5\textwidth}
			\centering
			\definecolor{qqqqff}{rgb}{0.,0.,1.}
			\begin{tikzpicture}[line cap=round,line join=round,>=triangle 45,x=1.0cm,y=1.0cm,scale=0.8]
			\clip(-4.589418683051358,-3.8254073116110523) rectangle (4.683985881573255,3.73721399712591);
			\draw [line width=1.6pt,color=qqqqff] (-3.381389745627314,-3.0148362354173153)-- (-2.235135560135566,-0.0240055331703113);
			\draw [line width=1.6pt,color=qqqqff] (-2.235135560135566,-0.0240055331703113)-- (2.,3.);
			\draw [line width=1.6pt,color=qqqqff] (2.,3.)-- (3.4712035376353554,1.8781074994199958);
			\draw [line width=1.6pt,color=qqqqff] (0.14532176812752454,-1.8392657308323692)-- (3.4712035376353554,1.8781074994199958);
			\draw [line width=1.6pt,color=qqqqff] (0.14532176812752454,-1.8392657308323692)-- (-3.381389745627314,-3.0148362354173153);
			\begin{scriptsize}
			\draw [fill=qqqqff] (3.4712035376353554,1.8781074994199958) circle (2.5pt);
			\draw [fill=qqqqff] (0.14532176812752454,-1.8392657308323692) circle (2.5pt);
			\draw [fill=qqqqff] (-3.381389745627314,-3.0148362354173153) circle (2.5pt);
			\draw [fill=qqqqff] (-2.235135560135566,-0.0240055331703113) circle (2.5pt);
			\draw [fill=qqqqff] (2.,3.) circle (2.5pt);
			\end{scriptsize}
			\end{tikzpicture}
		\end{minipage}%
		\begin{minipage}{.5\textwidth}
			\centering
			\definecolor{qqqqff}{rgb}{0.,0.,1.}
			\begin{tikzpicture}[line cap=round,line join=round,>=triangle 45,x=1.0cm,y=1.0cm,scale=1.5]
			\clip(-1.5651473732029693,-2.9106063502671127) rectangle (3.0622226760926092,0.8630937071756701);
			\draw [line width=1.6pt,color=qqqqff] (0.4226497308103742,0.7320508075688773)-- (-0.29708629022101113,0.5145685488949442);
			\draw [line width=1.6pt,color=qqqqff] (-0.29708629022101113,0.5145685488949442)-- (-0.6666666666666666,0.);
			\draw [line width=1.6pt,color=qqqqff] (-0.6666666666666666,0.)-- (-0.612004618869898,-1.0600230943494897);
			\draw [line width=1.6pt,color=qqqqff] (-0.612004618869898,-1.0600230943494897)-- (1.5773502691896257,-2.732050807568877);
			\draw [line width=1.6pt,color=qqqqff] (1.5773502691896257,-2.732050807568877)-- (2.,0.);
			\draw [line width=1.6pt,color=qqqqff] (2.,0.)-- (0.4226497308103742,0.7320508075688773);
			\begin{scriptsize}
			\draw [fill=qqqqff] (0.4226497308103742,0.7320508075688773) circle (1.5pt);
			\draw [fill=qqqqff] (-0.29708629022101113,0.5145685488949442) circle (1.5pt);
			\draw [fill=qqqqff] (-0.6666666666666666,0.) circle (1.5pt);
			\draw [fill=qqqqff] (-0.612004618869898,-1.0600230943494897) circle (1.5pt);
			\draw [fill=qqqqff] (1.5773502691896257,-2.732050807568877) circle (1.5pt);
			\draw [fill=qqqqff] (2.,0.) circle (1.5pt);
			\end{scriptsize}
			\end{tikzpicture}
		\end{minipage}
		\caption{An affinely regular pentagon (left) and a projectively regular hexagon (right).} \label{fig:regular}
	\end{figure}
	
	%[Finite subgroups of automorphism groups]
	
	%[Define regular polygons by intrinsic properties/quantities]
	
	\section{A projective analogue}
	In this section, we are concerned exclusively with projective geometry. The notions of natural constructions, regular polygons, moduli spaces of polygon, and so on, should be understood within the context of projective geometry. The words \textit{projective} and \textit{projectively} are sometimes omitted to avoid repetition.
	
	Given $4$ points $A, B, C, D$, we construct the points $P = AC \cap BD$, $Q = AB \cap CD$ and $H = BC \cap PQ$. Then we parametrize the line $PQ$ by using a projective transformation to map it to the projective line $\R P^1 = \R \cup \{\infty\}$ so that $P, H, Q$ correspond to the points $0, 1, \infty$, respectively. Let $H_{\lambda}$ denote the preimage of a point $\lambda$ in $\R P^1$. We write $H_{\lambda}(A,B,C,D)$ to denote this natural construction, which depends on the parameter $\lambda \in \R P^1$. 
	
	This construction can be extended cyclically to any polygon. More precisely, we can apply $H_{\lambda}$ to consecutive quadruples of vertices of a polygon and obtain a  natural polygon iteration, also denoted as $H_{\lambda}$.
	
	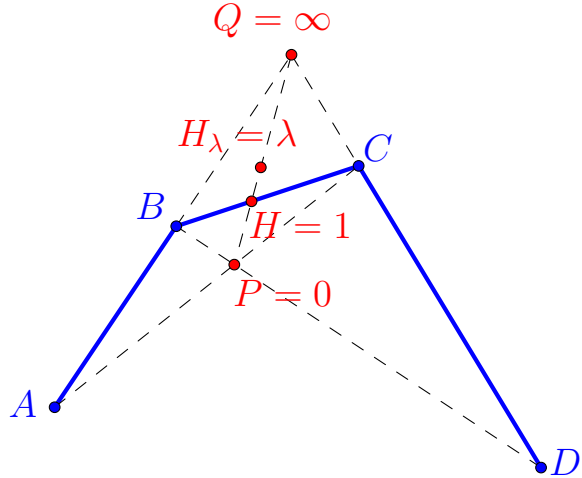
\begin{figure}[!ht]
		\centering
		\definecolor{ffqqqq}{rgb}{1.,0.,0.}
		\definecolor{qqqqff}{rgb}{0.,0.,1.}
		\begin{tikzpicture}[line cap=round,line join=round,>=triangle 45,x=1.0cm,y=1.0cm,scale=0.8]
		\clip(-5.88,-2.86) rectangle (5.72,6.6);
		\draw [line width=1.6pt,color=qqqqff] (-4.,-1.)-- (-2.,2.);
		\draw [line width=1.6pt,color=qqqqff] (-2.,2.)-- (1.,3.);
		\draw [line width=1.6pt,color=qqqqff] (1.,3.)-- (4.,-2.);
		\draw [dash pattern=on 5pt off 5pt] (-4.,-1.)-- (1.,3.);
		\draw [dash pattern=on 5pt off 5pt] (4.,-2.)-- (-2.,2.);
		\draw [dash pattern=on 5pt off 5pt] (-2.,2.)-- (-0.10526315789473684,4.842105263157895);
		\draw [dash pattern=on 5pt off 5pt] (-0.10526315789473684,4.842105263157895)-- (1.,3.);
		\draw [dash pattern=on 5pt off 5pt] (-0.10526315789473684,4.842105263157895)-- (-1.0454545454545454,1.3636363636363635);
		\begin{scriptsize}
		\draw [fill=qqqqff] (-4.,-1.) circle (2.5pt);
		\draw[color=qqqqff] (-4.52,-0.94) node {\large $A$};
		\draw [fill=qqqqff] (-2.,2.) circle (2.5pt);
		\draw[color=qqqqff] (-2.42,2.34) node {\large $B$};
		\draw [fill=qqqqff] (1.,3.) circle (2.5pt);
		\draw[color=qqqqff] (1.32,3.3) node {\large $C$};
		\draw [fill=qqqqff] (4.,-2.) circle (2.5pt);
		\draw[color=qqqqff] (4.4,-1.9) node {\large $D$};
		\draw [fill=ffqqqq] (-1.0454545454545454,1.3636363636363635) circle (2.5pt);
		\draw[color=ffqqqq] (-0.25,0.9) node {\large $P=0$};
		\draw [fill=ffqqqq] (-0.10526315789473684,4.842105263157895) circle (2.5pt);
		\draw[color=ffqqqq] (-0.42,5.4) node {\large $Q=\infty$};
		\draw [fill=ffqqqq] (-0.7619047619047619,2.4126984126984126) circle (2.5pt);
		\draw[color=ffqqqq] (0.04,2.02) node {\large $H=1$};
		\draw [fill=ffqqqq] (-0.6097560975609756,2.975609756097561) circle (2.5pt);
		\draw[color=ffqqqq] (-1.,3.5) node {\large $H_{\lambda}=\lambda$};
		\end{scriptsize}
		\end{tikzpicture}
		\caption{The construction $H_{\lambda}$.} \label{fig:Hlambda}
	\end{figure}
	
	Now we can state our main result.
	
	\begin{thm}[Main Theorem] \label{thm:main}
		Let $\phi = \frac{1+\sqrt{5}}{2}$ be the golden ratio. Given an arbitrary pentagon whose vertices are denoted as $A, B, C, D, E$, the pentagon with vertices 
		\begin{align*}
			A' &= H_{\phi}(B,C,D,E), \\ B' &= H_{\phi}(C,D,E,A), \\ C' &= H_{\phi}(D,E,A,B), \\ D' &= H_{\phi}(E,A,B,C), \\ E' &= H_{\phi}(A,B,C,D),
		\end{align*}
		 is projectively regular. In other words, the natural polygon iteration $H_{\phi}$ is immediately regularizing.
	\end{thm}
	
	\begin{figure}[!ht]
		\centering
		\definecolor{ffqqqq}{rgb}{1.,0.,0.}
		\definecolor{qqqqff}{rgb}{0.,0.,1.}
		\begin{tikzpicture}[line cap=round,line join=round,>=triangle 45,x=1.0cm,y=1.0cm,scale=0.8]
		\clip(-4.2709090909090905,-3.17454545454546) rectangle (8.838181818181825,5.898181818181816);
		\draw [dash pattern=on 5pt off 5pt] (-1.82,5.26)-- (0.58,-2.94);
		\draw [dash pattern=on 5pt off 5pt] (0.58,-2.94)-- (7.12,5.52);
		\draw [dash pattern=on 5pt off 5pt] (7.12,5.52)-- (-3.4,0.94);
		\draw [dash pattern=on 5pt off 5pt] (-3.4,0.94)-- (7.82,-0.78);
		\draw [dash pattern=on 5pt off 5pt] (7.82,-0.78)-- (-1.82,5.26);
		\draw [dash pattern=on 5pt off 5pt] (-0.42199371842141215,0.4834785379398243)-- (1.6003506295732268,3.1169587341678118);
		\draw [dash pattern=on 5pt off 5pt] (1.6003506295732268,3.1169587341678118)-- (2.8399559197513424,-0.016570782706979464);
		\draw [dash pattern=on 5pt off 5pt] (2.8399559197513424,-0.016570782706979464)-- (-0.8770864522497326,2.038378711853253);
		\draw [dash pattern=on 5pt off 5pt] (-0.8770864522497326,2.038378711853253)-- (4.067394924662966,1.5712172878667718);
		\draw [dash pattern=on 5pt off 5pt] (4.067394924662966,1.5712172878667718)-- (-0.42199371842141215,0.4834785379398243);
		\draw [dash pattern=on 5pt off 5pt] (0.2805821371982275,1.3983670454073764)-- (-3.4,0.94);
		\draw [dash pattern=on 5pt off 5pt] (0.6605089188012762,1.893104579280078)-- (-1.82,5.26);
		\draw [dash pattern=on 5pt off 5pt] (2.1397899107852827,1.7533400697163113)-- (7.12,5.52);
		\draw [dash pattern=on 5pt off 5pt] (2.3741337603716763,1.1609552391614828)-- (7.82,-0.78);
		\draw [dash pattern=on 5pt off 5pt] (1.2171015412799575,0.8806166720029581)-- (0.58,-2.94);
		\draw [line width=1.6pt,color=qqqqff] (-0.8770864522497326,2.038378711853253)-- (1.6003506295732268,3.1169587341678118);
		\draw [line width=1.6pt,color=qqqqff] (1.6003506295732268,3.1169587341678118)-- (4.067394924662966,1.5712172878667718);
		\draw [line width=1.6pt,color=qqqqff] (4.067394924662966,1.5712172878667718)-- (2.8399559197513424,-0.016570782706979464);
		\draw [line width=1.6pt,color=qqqqff] (2.8399559197513424,-0.016570782706979464)-- (-0.42199371842141215,0.4834785379398243);
		\draw [line width=1.6pt,color=qqqqff] (-0.42199371842141215,0.4834785379398243)-- (-0.8770864522497326,2.038378711853253);
		\draw [line width=1.6pt,color=ffqqqq] (-1.0284655316404216,1.2353427492778488)-- (-0.000008902150378962726,2.789652223115915);
		\draw [line width=1.6pt,color=ffqqqq] (-0.000008902150378962726,2.789652223115915)-- (3.238682521850525,2.584460583035625);
		\draw [line width=1.6pt,color=ffqqqq] (3.238682521850525,2.584460583035625)-- (3.9293495660130575,0.6066624722880338);
		\draw [line width=1.6pt,color=ffqqqq] (-1.0284655316404216,1.2353427492778488)-- (1.0623413987476287,-0.047460160979403775);
		\draw [line width=1.6pt,color=ffqqqq] (3.9293495660130575,0.6066624722880338)-- (1.0623413987476287,-0.047460160979403775);
		\begin{scriptsize}
		\draw [fill=qqqqff] (-0.42199371842141215,0.4834785379398243) circle (2.5pt);
		\draw[color=qqqqff] (-0.8709090909090889,0.24363636363635915) node {\large $A$};
		\draw [fill=qqqqff] (1.6003506295732268,3.1169587341678118) circle (2.5pt);
		\draw [fill=qqqqff] (2.8399559197513424,-0.016570782706979464) circle (2.5pt);
		\draw [fill=qqqqff] (-0.8770864522497326,2.038378711853253) circle (2.5pt);
		\draw [fill=qqqqff] (4.067394924662966,1.5712172878667718) circle (2.5pt);
		\draw [fill=ffqqqq] (3.238682521850525,2.584460583035625) circle (2.5pt);
		\draw[color=ffqqqq] (3.1290909090909134,3.0618181818181784) node {\large $A'$};
		\draw [fill=ffqqqq] (-0.000008902150378962726,2.789652223115915) circle (2.5pt);
		\draw [fill=ffqqqq] (-1.0284655316404216,1.2353427492778488) circle (2.5pt);
		\draw [fill=ffqqqq] (1.0623413987476287,-0.047460160979403775) circle (2.5pt);
		\draw [fill=ffqqqq] (3.9293495660130575,0.6066624722880338) circle (2.5pt);
		\end{scriptsize}
		\end{tikzpicture}
		\caption{The map $H_{\phi}$ on a pentagon.} \label{fig:Hphi}
	\end{figure}
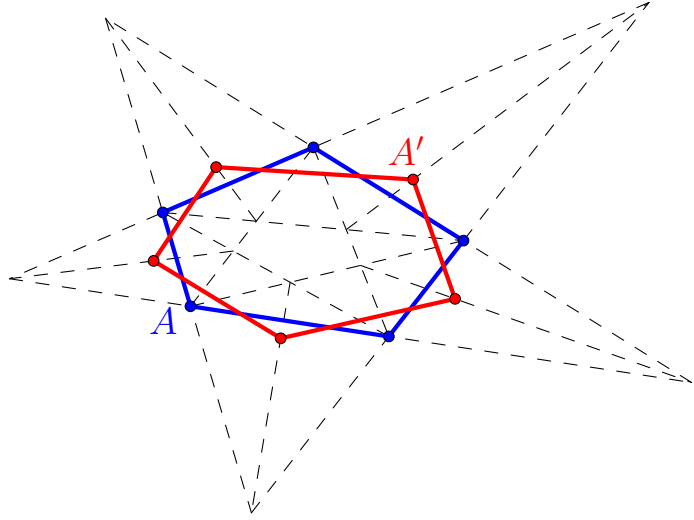
	
	Note that the Main Theorem, as well as Napoleon's and Varignon's theorems, do not require the involved polygons to be convex or simple. We also note that the points $A, B, H_{\phi}(A,B,C,D), C, D$ form a projectively regular pentagon.
	
	One way to prove this theorem is that we first symbolically compute the points $A',B',C',D',E'$ based on the coordinates of $A,B,C,D,E$ and then check that the pentagon $A'B'C'D'E'$ is indeed projectively regular. This can be done with any computer algebra system such as Sage, Maple and Mathematica. A purely geometric proof would be desirable but is not known to the author.
	
	The Main Theorem was discovered in our study of the polygon iterations $H_{\lambda}$ on pentagons. Empirically, for $\lambda \neq 0,\infty$, it is observed that $H_{\lambda}$ is regularizing for almost all starting pentagons. Therefore, we can try to define the Julia set of $H_{\lambda}$ to be the subset of the moduli space of pentagons that will not converge to the equivalence classes of regular polygons. Here are the pictures of the Julia set of $H_{\lambda}$ for some $\lambda$ in $(0,\infty)$.
	
	\begin{figure}[!ht]
		\centering
		\begin{minipage}{.3\textwidth}
			\centering
			\includegraphics[width=4cm]{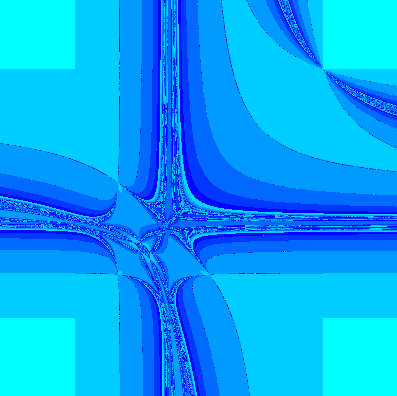}
		\end{minipage}%
		\begin{minipage}{.3\textwidth}
			\centering
			\includegraphics[width=4cm]{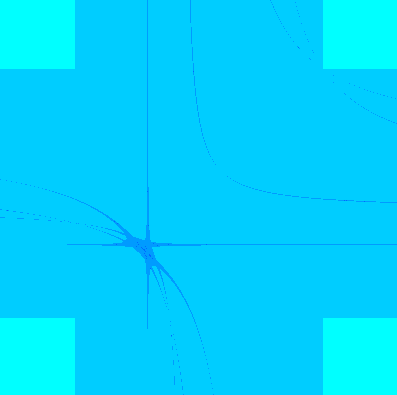}
		\end{minipage}%
		\begin{minipage}{.3\textwidth}
			\centering
			\includegraphics[width=4cm]{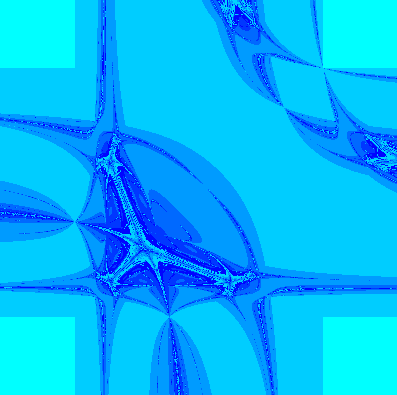}
		\end{minipage}%
		\caption{The Julia sets of $H_{\lambda}$ with $\lambda = 0.2$, $1.4$ ($\approx \phi$) and $7.0$ in an appropriate coordinate system of the ($2$-dimensional) moduli space of pentagons.}
	\end{figure}
	
	Observe that the Julia set seems to disappear when $\lambda$ tends to the golden ratio $\phi = \frac{1+\sqrt{5}}{2}$. This made us notice the distinctiveness of this particular parameter and led to the finding of the Main Theorem. The polygon iterations $H_{\lambda}$ are themselves fascinating objects of study and are explored in \cite{le2016family}. The best-known members of this family are perhaps the pentagram map $H_0$ and its inverse $H_{\infty}$, which are discrete integrable systems \cite{ovsienko2010pentagram, beffa2013generalizations}. Another prominent member is the projective heat map $H_1$, which has been studied extensively for the case of pentagons in \cite{schwartz2016heat} and was the original motivation for the maps $H_{\lambda}$.
	
	There is a related result.
	
	\begin{thm}
		Let $\psi = \frac{1-\sqrt{5}}{2} = -\frac{1}{\phi}$. Given an arbitrary pentagon $ABCDE$, the pentagon with vertices 
		\begin{align*}
		A' &= H_{\psi}(B,C,D,E), \\ B' &= H_{\psi}(C,D,E,A), \\ C' &= H_{\psi}(D,E,A,B), \\ D' &= H_{\psi}(E,A,B,C), \\ E' &= H_{\psi}(A,B,C,D),
		\end{align*}
		is projectively star-regular, i.e. it is projectively equivalent to a star-regular pentagon.
	\end{thm}
	
	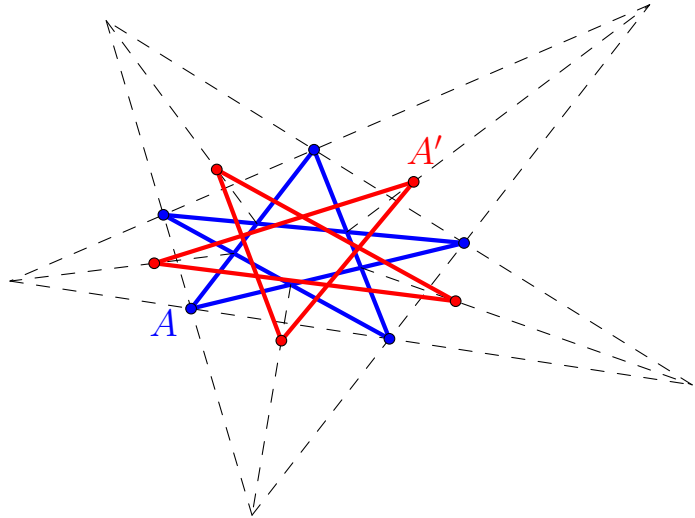
\begin{figure}[!ht]
		\centering
		\definecolor{ffqqqq}{rgb}{1.,0.,0.}
		\definecolor{qqqqff}{rgb}{0.,0.,1.}
		\begin{tikzpicture}[line cap=round,line join=round,>=triangle 45,x=1.0cm,y=1.0cm,scale=0.8]
		\clip(-4.512677617093572,-3.331734915036172) rectangle (8.731176517505222,5.83426121972776);
		\draw [dash pattern=on 5pt off 5pt] (-1.82,5.26)-- (0.58,-2.94);
		\draw [dash pattern=on 5pt off 5pt] (0.58,-2.94)-- (7.12,5.52);
		\draw [dash pattern=on 5pt off 5pt] (7.12,5.52)-- (-3.4,0.94);
		\draw [dash pattern=on 5pt off 5pt] (-3.4,0.94)-- (7.82,-0.78);
		\draw [dash pattern=on 5pt off 5pt] (7.82,-0.78)-- (-1.82,5.26);
		\draw [line width=1.6pt,color=qqqqff] (-0.42199371842141215,0.4834785379398243)-- (1.6003506295732268,3.1169587341678118);
		\draw [line width=1.6pt,color=qqqqff] (1.6003506295732268,3.1169587341678118)-- (2.8399559197513424,-0.016570782706979464);
		\draw [line width=1.6pt,color=qqqqff] (2.8399559197513424,-0.016570782706979464)-- (-0.8770864522497326,2.038378711853253);
		\draw [line width=1.6pt,color=qqqqff] (-0.8770864522497326,2.038378711853253)-- (4.067394924662966,1.5712172878667718);
		\draw [line width=1.6pt,color=qqqqff] (4.067394924662966,1.5712172878667718)-- (-0.42199371842141215,0.4834785379398243);
		\draw [dash pattern=on 5pt off 5pt] (0.2805821371982275,1.3983670454073764)-- (-3.4,0.94);
		\draw [dash pattern=on 5pt off 5pt] (0.6605089188012762,1.893104579280078)-- (-1.82,5.26);
		\draw [dash pattern=on 5pt off 5pt] (2.1397899107852827,1.7533400697163113)-- (7.12,5.52);
		\draw [dash pattern=on 5pt off 5pt] (2.3741337603716763,1.1609552391614828)-- (7.82,-0.78);
		\draw [dash pattern=on 5pt off 5pt] (1.2171015412799575,0.8806166720029581)-- (0.58,-2.94);
		\draw [line width=1.6pt,color=ffqqqq] (3.238682521850525,2.584460583035626)-- (-1.0284655316404214,1.2353427492778484);
		\draw [line width=1.6pt,color=ffqqqq] (-1.0284655316404214,1.2353427492778484)-- (3.9293495660130566,0.6066624722880336);
		\draw [line width=1.6pt,color=ffqqqq] (3.9293495660130566,0.6066624722880336)-- (-0.000008902150378962726,2.789652223115915);
		\draw [line width=1.6pt,color=ffqqqq] (-0.000008902150378962726,2.789652223115915)-- (1.062341398747629,-0.04746016097940388);
		\draw [line width=1.6pt,color=ffqqqq] (1.062341398747629,-0.04746016097940388)-- (3.238682521850525,2.584460583035626);
		\begin{scriptsize}
		\draw [fill=qqqqff] (-0.42199371842141215,0.4834785379398243) circle (2.5pt);
		\draw[color=qqqqff] (-0.8756691307543742,0.25016738211606493) node {\large $A$};
		\draw [fill=qqqqff] (1.6003506295732268,3.1169587341678118) circle (2.5pt);
		\draw [fill=qqqqff] (2.8399559197513424,-0.016570782706979464) circle (2.5pt);
		\draw [fill=qqqqff] (-0.8770864522497326,2.038378711853253) circle (2.5pt);
		\draw [fill=qqqqff] (4.067394924662966,1.5712172878667718) circle (2.5pt);
		\draw [fill=ffqqqq] (3.238682521850525,2.584460583035626) circle (2.5pt);
		\draw[color=ffqqqq] (3.4226136258283133,3.1136331928190113) node {\large $A'$};
		\draw [fill=ffqqqq] (-1.0284655316404214,1.2353427492778484) circle (2.5pt);
		\draw [fill=ffqqqq] (3.9293495660130566,0.6066624722880336) circle (2.5pt);
		\draw [fill=ffqqqq] (-0.000008902150378962726,2.789652223115915) circle (2.5pt);
		\draw [fill=ffqqqq] (1.062341398747629,-0.04746016097940388) circle (2.5pt);
		\end{scriptsize}
		\end{tikzpicture}
		\caption{The map $H_{\psi} = H_{-1/\phi}$ on a pentagon.} \label{fig:Hpsi}
	\end{figure}
	
	Yes, you see it right -- the sets $\{A,B,C,D,E\}$ and $\{A',B',C',D',E'\}$ are the same in both Figure \ref{fig:Hphi} and Figure \ref{fig:Hpsi}! This is not a coincidence but the reflection of a symmetry within the family $H_{\lambda}$. Specifically, $H_{\lambda}$ and $H_{-1/\lambda}$ are conjugate modulo a relabeling of the vertices that transforms a convex pentagon to a star-convex one. This fact is not hard to show, but see \cite{le2016family} for more details.
	
	\section{Some further questions}
	Here are some interesting and related questions.
	
	\begin{que}
		Are there higher-dimensional analogues of Napoleon's theorem, Varignon's theorem and the Main Theorem? 
	\end{que}
	
	\begin{que}
		Are there analogues in other Klein geometries?
	\end{que}
	
	\begin{que}
		Are there analogues of Varignon's theorem and the Main Theorem for other $n$-gons in the spirit of the Petr-Neumann-Douglas theorem \cite{Gruenbaum2001}?
	\end{que}
	
	\section{Acknowledgment}
	The author is deeply grateful to his advisor, Richard E. Schwartz, for introducing him to these fascinating maps and for suggesting the connection between the Main Theorem and Varignon's theorem. He would also like to thank Sergei Tabachnikov and Sinai Robins for many insightful conversations.


\begin{thebibliography}{}
		
		\bibitem[\protect\astroncite{Gr{\"u}nbaum}{2001}]{Gruenbaum2001}
		Gr{\"u}nbaum, B. (2001).
		\newblock A relative of “Napoleon's theorem”.
		\newblock {\em Geombinatorics}, 10:116--121.
		
		\bibitem[\protect\astroncite{Klein}{1893}]{Klein1893}
		Klein, F. (1893).
		\newblock Vergleichende betrachtungen {\"u}ber neuere geometrische forschungen.
		\newblock {\em Mathematische Annalen}, 43(1):63--100.
		
		\bibitem[\protect\astroncite{Le}{2016}]{le2016family}
		Le, Q.-N. (2016).
		\newblock A family of projectively natural polygon iterations.
		\newblock {\em arXiv preprint arXiv:1602.02699}.
		
		\bibitem[\protect\astroncite{Mar{\'\i}-Beffa}{2013}]{beffa2013generalizations}
		Mar{\'\i}-Beffa, G. (2013).
		\newblock On generalizations of the pentagram map: discretizations of agd
		flows.
		\newblock {\em Journal of nonlinear science}, 23(2):303--334.
		
		\bibitem[\protect\astroncite{Ovsienko et~al.}{2010}]{ovsienko2010pentagram}
		Ovsienko, V., Schwartz, R., and Tabachnikov, S. (2010).
		\newblock The pentagram map: a discrete integrable system.
		\newblock {\em Communications in Mathematical Physics}, 299(2):409--446.
		
		\bibitem[\protect\astroncite{Schwartz}{2016}]{schwartz2016heat}
		Schwartz, R. (2016).
		\newblock The projective heat map.
		\newblock {\em Preprint, available at } \url{http://www.math.brown.edu/~res/Papers/heat2.pdf}.
		
		\bibitem[\protect\astroncite{Tabachnikov}{1999}]{Tabachnikov1999}
		Tabachnikov, S. (1999).
		\newblock Fagnano orbits of polygonal dual billiards.
		\newblock {\em Geometriae Dedicata}, 77(3):279--286.
		
		\bibitem[\protect\astroncite{Terras}{1999}]{Terras1999}
		Terras, A. (1999).
		\newblock {\em Fourier analysis on finite groups and applications}.
		\newblock Number~43. Cambridge University Press.
		
		\bibitem[\protect\astroncite{Troubetzkoy}{2009}]{Troubetzkoy2009}
		Troubetzkoy, S. (2009).
		\newblock Dual billiards, Fagnano orbits, and regular polygons.
		\newblock {\em American Mathematical Monthly}, 116(3):251--260.
		
	\end{thebibliography}
\end{document}